\documentclass[]{siamart250211}
\usepackage[utf8]{inputenc}
\usepackage[english]{babel}
\AddToHook{cmd/appendix/before}{%
    \crefalias{section}{appendix}%
    \crefalias{subsection}{appendix}
}
\usepackage{amsfonts}
\usepackage{algorithmic}
\renewcommand{\algorithmicrequire}{\textbf{input:}}
\renewcommand{\algorithmicensure}{\textbf{output:}}
\algsetup{
  linenodelimiter=\ ,
  linenosize=\normalsize
}
\Crefname{ALC@unique}{Line}{Lines}
\usepackage{mathtools}
\numberwithin{equation}{section}
\usepackage{xfrac}
\usepackage[matrix,arrow]{xy}
\usepackage{svg}
\usepackage{xcolor}

\newcommand{\calC}{\mathcal{C}}

\newcommand{\bbN}{\mathbb{N}}
\newcommand{\bbR}{\mathbb{R}}
\newcommand{\bbZ}{\mathbb{Z}}
\newcommand{\frakg}{\mathfrak{g}}
\newcommand{\frakgl}{\mathfrak{gl}}
\newcommand{\frakso}{\mathfrak{so}}
\newcommand{\diff}{\mathrm{d}}

\newcommand{\dd}[2][{}]{\frac{\!\diff#1}{\!\diff#2}}
\newcommand{\ddt}{\dd{t}}
\newcommand{\pp}[2][{}]{\frac{\partial#1}{\partial#2}}

\makeatletter
\newcommand{\ros}{\mathop{\operator@font ros}\nolimits}
\newcommand{\dros}{\mathop{\operator@font dros}\nolimits}
\newcommand{\cay}{\mathop{\operator@font cay}\nolimits}
\newcommand{\dcay}{\mathop{\operator@font dcay}\nolimits}
\newcommand{\dexp}{\mathop{\operator@font dexp}\nolimits}
\newcommand{\dlog}{\mathop{\operator@font dlog}\nolimits}
\renewcommand{\skew}{\mathop{\operator@font skew}\nolimits}
\newcommand{\unskew}{\mathop{\operator@font unskew}\nolimits}
\newcommand{\dskew}{\mathop{\operator@font dskew}\nolimits}
\newcommand{\dunskew}{\mathop{\operator@font dunskew}\nolimits}
\newcommand{\argmin}{\mathop{\operator@font argmin}}
\newcommand{\tr}{\mathop{\operator@font tr}\nolimits}
\newcommand{\la}{\mathop{\operator@font L}\nolimits}
\newcommand{\ra}{\mathop{\operator@font R}\nolimits}
\newcommand{\lra}{\mathop{\operator@font A}\nolimits}
\newcommand{\ad}{\mathop{\operator@font ad}\nolimits}
\newcommand{\Ad}{\mathop{\operator@font Ad}\nolimits}
\newcommand{\rtr}{\tau}
\newcommand{\rct}{{\tau^{-1}}}
\newcommand{\drtr}{{\operator@font d}\rtr}
\newcommand{\drct}{{\operator@font d}\rct}
\newcommand{\Tan}{{\operator@font T}}
\makeatother

\newcommand{\norm}[1]{\|#1\|}
\newcommand\gopt{g^\star}
\newcommand{\coloneq}{\vcentcolon \mathrel {\mkern -1.2mu}\mathrel {=}}
\newcommand{\eqcolon}{\mathrel {=}\mathrel {\mkern -1.2mu}\vcentcolon}

\newcommand{\To}{\longrightarrow}
\newcommand{\assign}{\leftarrow}
\newcommand{\quand}{\quad\text{and}\quad}
\newcommand{\qquand}{\qquad\text{and}\qquad}

\title{%
  Momentum-based gradient descent methods\texorpdfstring{\\}{ }for Lie groups%
  \thanks{%
    Submitted to the editors on \today.%
    \funding{%
      CMC and DMdD acknowledge financial support from the Spanish Ministry of Science and Innovation under grants PID2022-137909NB-C21, TED2021-129455B-I00, and CEX2023-001347-S funded by MCIN/AEI\-/10.13039\-/\-501100011033. DMdD acknowledges financial support from BBVA Foundation via the project ``Mathematical optimization for a more efficient, safer and decarbonized maritime transport''.}}%
}
\author{%
  Cédric M. Campos%
  \thanks{%
    Departament de Matemàtiques, %
    Universitat de València, %
    Carrer Dr. Moliner 50, 46100 Burjassot, Spain %
    (\email{cedric.martinez@uv.es}, \url{https://cmcampos.xyz}). %
    \emph{On leave from} %
    Departamento de Matemática Aplicada, Ciencia e Ingeniería %
    de los Materiales y Tec\-no\-lo\-gía Electrónica, %
    Universidad Rey Juan Carlos, %
    Calle Tulipán s/n, 28933 Móstoles, Spain.
  }%
  \and%
  David Martín de Diego%
  \thanks{%
    Instituto de Ciencias Matemáticas (CSIC-UAM-UC3M-UCM), %
    Calle Nicolás Cabrera 13-15, 28049 Madrid, Spain %
    (\email{david.martin@icmat.es}).%
  }%
  \and%
  Jose Torrente-Teruel%
  \thanks{%
    Departamento de Matemáticas, %
    Universidad de Córdoba, %
    Edificio Albert Einstein, Campus de Rabanales, %
    14071 Córdoba, Spain %
    (\email{jtorrente@uco.es}).
  }%
}
\headers{
  Momentum-based gradient descent methods for Lie groups%
}{%
  C.~M.~Campos, D.~Martín de Diego, J.~Torrente%
}

\begin{document}
\maketitle

\begin{abstract}%
  Polyak’s Heavy Ball (PHB) \cite{Po64}, also known as Classical Momentum, and Nesterov’s Accelerated Gradient (NAG) \cite{Ne83} are well-established momentum-descent methods for optimization. Although the latter generally outperforms the former, primarily, generalizations of PHB-like methods to nonlinear spaces have not been sufficiently explored in the literature. In this paper, we propose a generalization of NAG-like methods for Lie group optimization. This generalization is based on the variational one-to-one correspondence between classical and accelerated momentum methods \cite{CMM23}. We provide numerical experiments for chosen retractions on the group of rotations based on the Frobenius norm and the Rosenbrock function to demonstrate the effectiveness of our proposed methods, and that align with results of the Euclidean case, that is, a faster convergence rate for NAG.
\end{abstract}



\section{Introduction}
\label{sec:intro}
A fundamental step of many of the recent advances in machine learning and data analysis consists of the minimization of a loss function. This loss function allows us to evaluate, for instance, how well the machine learning algorithm models the featured data set. Due to the typically large size of data, low-cost optimization techniques such as the gradient descent (GD) method are more convenient than methods that require the computation of second-order derivatives, like Newton's method. Therefore, it is useful to accelerate gradient descent without increasing computational cost \cite{Nesterov2004}. Polyak \cite{Po64} introduced Classical Momentum (CM), also known as Polyak's Heavy Ball (PHB), as a technique to accelerate gradient descent by taking into account previous gradients in the update rule at each iteration of the method. Later, Nesterov \cite{Ne83} found Nesterov's Accelerated Gradient (NAG) method as an alternative optimization technique with an optimal convergence rate for the class of convex loss functions with Lipschitz gradient. All of these families of accelerated optimization methods have become popular in the machine learning community.

Given a convex function \(f\in\calC^2(\bbR^d,\bbR)\) and the corresponding minimization problem
\[
\argmin_{x\in\bbR^d} f(x)\,,
\]
observe that the different convergence behavior of GD and accelerated optimization is retained in the continuous limit of these methods \cite{SuBoCa16}:
\begin{align*}
  \dot{x}+\nabla f(x)&=0 \tag{GD} \,,\\
\ddot{x}+\frac{3}{t}\dot{x}+\nabla f(x)&=0 \tag{PHB/NAG}\,.
\end{align*}
GD is modeled by a first-order differential equation, while the continuous limit of accelerated methods such as PHB and NAG consists of a second-order differential equation (SODE). This SODE can be recovered from a variational principle as the Euler-Lagrange equations for the time-dependent Lagrangian \cite{WiWiJo16}
\[
  L(t,x,\dot{x}) = t^3\left(\frac12||\dot{x}||^2-f(x)\right) \, .
\]
However, a force must be included to obtain the NAG method, hence modifying the SODE \cite{CMM23}. The simulation of Lagrangian or Hamiltonian systems has made it possible to introduce discrete variational \cite{marsden-west} and symplectic methods \cite{serna,hairer,blanes} as a sub-product of the classical accelerated optimization methods. In particular, Campos et al. \cite{CMM23} introduced variational integrators which allowed to generalize PHB and NAG, deriving two families of optimization methods in one-to-one correspondence. However, since the systems considered are explicitly time-dependent, the preservation of symplecticity occurs solely on the fibers.

In the majority of machine learning applications, the function to be optimized is modeled on a Euclidean space but other cases are also of considerable interest (see \cite{DuruLeok2022,DuruLeok2023,LeeTaoLeok,Sharma} and references therein). Particularly, in this paper we study optimization problems where the objective function is defined on a Lie group \cite{AbMaSeBookRetraction} as in signal or image processing, independent component analysis (ICA), learning robotic systems etc (see \cite{Bernardini2021,TaOh20,DuruLeok2023} and references therein). Such problems are usually tackled using similar techniques as in the standard Euclidean case, using, for instance, a constrained optimization procedure or an appropriate parametrization to transform them into unconstrained problems. Such algorithms are characterized by a reduced convergence due to the lack of a geometric framework. In this paper, we adopt an intrinsic point of view, constructing the accelerated methods on Lie groups using its inherent geometry. In addition, the left/right trivialization is used as a fundamental tool in order to simplify and obtain more efficient methods, in contrast to general differential manifold structures. In arbitrary manifolds it is necessary to use more involved techniques, as for instance, to equip the manifold with a Riemannian metric and define a retraction map from it or using projections from an euclidean space (see \cite{AbMaSeBookRetraction} for more details). However, defining such general methods on manifolds is complicated, and in the case of Lie groups we can use the left/right trivializations to simplify the geometry to a vector space (the Lie algebra). In particular, in this work we introduce PHB-type methods on Lie groups without relying on an extended Lagrangian formalism, as used in \cite{LeeTaoLeok}. Furthermore, we derive a NAG-type extension to Lie groups by incorporating appropriate external forces. According to our derivation, and in contrast with some interpretations in the literature, existing momentum-based methods on Lie groups are more accurately classified within the PHB family (see, for instance, \cite{TaOh20}).

The paper is organized as follows. In \Cref{sec:the methods}, we introduce the notation to be used in the following and give schematically the algorithms developed in this work. In fact, PHB and NAG methods in Lie groups can be computed using \cref{alg:gd:G}. \Cref{sec:devel} is devoted to the derivation of both method families, Eqs. \cref{eq:mm:G}, using a discrete variational perspective from a forced discrete Lagrangian system on a Lie group. We also give an alternative derivation from a Hamilton-Pontryagin variational principle. In the remaining two sections are devoted to exemplify the methods and test the computational performance of the optimization techniques with respect to the Gradient Descent method. Several objective functions are defined, and explicit solvers for these (\textit{a priori} implicit) methods are presented in \Cref{sec:examples}. They involve two important retraction maps: the exponential map and the Cayley transform. Then, \Cref{sec:experiments} provides numerical simulations to the test functions. The algorithms introduced here are generally shown to be improvements over Gradient Descent, except for discrepancies in some cases. Conclusions and overall discussion can be found in \Cref{sec:conclusions}.
Finally, to make the paper self-contained, we include several appendices at the end, containing the necessary technical results and background theory on Lie groups and Discrete Geometric Mechanics used throughout the work. Most of these results can be found scattered across the literature, sometimes with divergent notations. For this reason, we present them here with a unified notation and complete proofs for those in \crefrange{sec:cayley}{sec:skew}. For further details, the reader is referred to \cite{AbMaSeBookRetraction,Bou1,IserlesEtAl00,Le13-SM,Le18-RM,marsden-west}.



\section{The methods}
\label{sec:the methods}
The sole purpose of this section is to present in a concise manner the family of methods whose derivation is developed in the next section. This allows the reader to immediately recognize the analogy with Euclidean PHB/NAG methods and facilitate their direct implementation. Before proceeding, we briefly summarize the notation. Some definitions are either assumed (see \cite{Le13-SM,Le18-RM}) or introduced later.

\subsection{Notation}
\label{sec:notation}

\begin{itemize}
\item \(G\) denotes a Lie group, the associated Lie algebra is  \(\frakg=\Tan_eG\), and \(\frakg^*\) its dual.
\item \(\la_g\) and \(\ra_h\) are the left and right actions of the group, \(\la_g(h)=gh=\ra_h(g)\). Their tangent maps at the identity, \(\Tan_e\la_g\) and \(\Tan_e\ra_h\), are still denoted \(\la_g\) and \(\ra_h\). In addition, the adjoint map is \(\Ad(g)=T_e(\la_g\circ\ra_{g^{-1}})\).
\item Given a real-valued function \(\phi\colon G\to\bbR\), \(\diff\phi\colon TG\to\bbR\) is the differential of \(\phi\), a 1-form over \(G\).
\item \((\cdot)^*\) denotes the pullback.
\item We consider an inner product \(\langle\,\cdot\!\;,\cdot\!\;\rangle\) on \(\frakg\), for which \((\cdot)^\flat\colon\frakg\to\frakg^*\) and \((\cdot)^\sharp\colon\frakg^*\to\frakg\) denote the musical operators, and \((\cdot)^t\) the transposition of linear maps.
\item \(\nabla\phi\) is the right-trivialized gradient, \(\nabla\phi(g) \coloneqq \left( \ra_g^*\diff\phi(g) \right)^\sharp\).
\item \(\rtr\colon\frakg\to G\) is a retraction map, and \(\diff\rtr_\xi\colon\frakg\to\frakg\), for \(\xi\in\frakg\), denotes its right-trivialized tangent (see \Cref{sec:lie.groups}).
\item \(\Delta\) is the forward difference operator. For vectors (and covectors), it is the standard operator, \emph{e.g.} \(\Delta[\omega_0] = \omega_1-\omega_0\), either in \(\frakg\) or in \(\frakg^*\). For group elements, it gives the right-transition, \(\Delta w_0=w_0^{-1}w_1\) in \(G\), an ``arrow'' pointing from \(w_0\) to \(w_1\) when acting on the right of \(w_0\): \(\ra_{\Delta w_0}(w_0)=w_0\Delta w_0=w_1\).
\end{itemize}

\subsection{Momentum-Descent Methods for Lie groups}
\label{sec:method}

Given a Lie group \(G\), let \(\phi\colon D\subseteq G\to\bbR\) denote a real-valued \(\calC^1\)-function defined on a path-connected open subset \(D\subseteq G\). Assume that \(\phi\) possesses a single local minimum in \(D\),
\begin{equation*}
  \textstyle
  \gopt = \argmin_{g\in D} \phi(g)\,.
\end{equation*}
To seek for \(\gopt\), we propose a family of twin methods inspired by the one-to-one correspondence between PHB and NAG methods \cite{CMM23}. In fact, they are equivalent to ``regular'' PHB and NAG when \(G=\bbR^n\). For further details, see \Cref{sec:devel}. This correspondence allows for the compilation of both in a single algorithm, \cref{alg:mm:G}, with a Boolean input or hyperparameter to set the family of choice: \(\epsilon=0\), PHB-like method; \(\epsilon=1\), NAG-like method. A further hyperparameter is the strategy, a sequences of couples of coefficients, \((\mu_k,\eta_k)\), \(k \in \bbN_0\coloneqq\bbN\cup\{0\}\). \(\mu_k\) is usually referred to as the momentum coefficient and \(\eta_k\) as the learning rate. There are more general strategies where \(\mu\) and \(\eta\) depend on \(\nabla\phi\) or the past trajectory, as the original method by Nesterov \cite{Ne83}, but such strategies are out of the scope of the present work. See \Cref{sec:devel} for the choice of strategy.

\begin{algorithm}
  \caption{Momentum-based gradient descent method for Lie groups. Minimizes \(\phi\) from the initial guess \(g_0\) with strategy \((\eta,\mu)\). Set \(\varepsilon=0\) for PHB, or \(\varepsilon=1\) for NAG.}
  \label{alg:mm:G}
  \centering
  \begin{minipage}{.75\linewidth}
    \begin{algorithmic}[1]
      \ttfamily\bfseries%
      \makeatletter
      \addtocounter{ALC@unique}{-1}
      \addtocounter{ALC@line}{-1}
      \makeatother
      \STATE\label{alg:mm:G:input}%
      \algorithmicrequire\ %
      \(\nabla\phi\colon G\to\frakg\), \(g_0\in G\); %
      \(\eta,\mu\colon\bbN_0\to\bbR\), \(\varepsilon\in\{0,1\}\)%
      \STATE\label{alg:mm:G:ini}%
      \(g_1\assign g_0\), \(x_0\assign0\), \(x_1\assign0\), %
      \(y_1 \assign -\eta_0\nabla\phi(g_0)\), \(z_1\assign \varepsilon y_1\)%
      \FOR{\(k = 1\) \TO \(N-1\)}\label{alg:mm:G:loop}%
      \STATE\label{alg:mm:G:step:gd}%
      \(y_{k+1} \assign x_k - \eta_k\nabla\phi(g_k)\)%
      \STATE\label{alg:mm:G:z}%
      \(z_{k+1} \assign (1-\varepsilon)x_k + \varepsilon y_{k+1}\)%
      \STATE\label{alg:mm:G:step:mom}%
      \(x_{k+1} \assign y_{k+1} + \mu_k\Delta z_k\)%
      \STATE\label{alg:mm:G:solve}%
      \(g_{k+1} \assign g_k\tau(\xi_k)\) such that %
      \(\xi_k = \diff\rtr_{\xi_k}^t\left( \Ad_{g_k}^t\Delta x_k \right)\)%
      \ENDFOR\label{alg:mm:G:end}%
      \STATE \algorithmicensure\ \(g_N\)
    \end{algorithmic}
  \end{minipage}
\end{algorithm}

The inputs are specified in \cref{alg:mm:G:input}, namely, \(\nabla\phi\), the right-trivialized gradient of the objective function, and \(g_0\), an initial guess for the minimizer. In \cref{alg:mm:G:ini}, the search direction is initialized to a safe value (stationary start), and several variables are set according to \cref{eq:mm:G:inival}. Beginning at \cref{alg:mm:G:loop} with \(k=1\), a gradient descent step is performed in \cref{alg:mm:G:step:gd}, followed by a momentum step in \cref{alg:mm:G:step:mom}. This yields a new momentum \(\Delta x_1\), which is then used together with \(g_1\) in \cref{alg:mm:G:solve} to compute a new approximation \(g_2\) of \(\gopt\) via the reconstruction equation \cref{eq:DEL:G:reconstruction}. This process is iterated through \crefrange{alg:mm:G:loop}{alg:mm:G:end}, following the dynamical equation \cref{eq:DEL:G} in the form of \cref{eq:mm:G}. The final iterate, \(g_N\), is then returned.

The variable of interest is \(g\); in fact, the sequence \(\{g_k\}\) is a trajectory of group elements converging toward \(\gopt\). The variables \(x\) and \(y\) are auxiliary elements in \(\frakg\) that carry part of the dynamics. The variable \(z\), introduced in \cref{alg:mm:G:z}, is an additional auxiliary variable in \(\frakg\) used to select the method family via the Boolean hyperparameter \(\varepsilon\). In a final implementation, according to the chosen family, either \(x\) or \(y\) should replace \(z\) in \cref{alg:mm:G:solve}. It is then readily seen that the steps in \crefrange{alg:mm:G:step:gd}{alg:mm:G:step:mom} resemble those of PHB/NAG methods.

The computational load is concentrated in the gradient evaluation, \cref{alg:mm:G:step:gd}, one per iteration, and in the reconstruction step, \cref{alg:mm:G:solve}. Although it is implicit in general, it can be rendered explicit in some cases. For instance, when \(G\) is the Euclidean space \(\bbR^n\), then \(g_k=x_k\) and \cref{alg:mm:G:solve} reduces to the tautological relation
\[ x_{k+1} = x_k + \Delta x_k \,. \]
And more notably, when \(G\) is the group of rotations \(SO(3)\) and \(\tau\) is the matrix exponential, then \(g_k=R_k\in SO(3)\) and the aforementioned equation reads as in Equation \cref{eq:solver:exp}, that is,
\[ R_{k+1} = \exp(\Delta x_k)R_k \,, \]
where \(\exp(\Delta x_k)\) is the exponential of a skewsymmetric matrix.

Finally, note that, for a strategy with zero momentum, \(\mu\equiv0\), we recover gradient descent for Lie groups, \cref{alg:gd:G}, which could be further simplified, but is left as is for easier comparison with \cref{alg:mm:G}.

\begin{algorithm}
  \caption{Gradient Descent for Lie groups. Minimizes \(\phi\) from the initial guess \(g_0\) with strategy \(\eta\).}
  \label{alg:gd:G}
  \centering
  \begin{minipage}{.68\linewidth}
    \begin{algorithmic}[1]
      \ttfamily\bfseries%
      \makeatletter
      \addtocounter{ALC@unique}{-1}
      \addtocounter{ALC@line}{-1}
      \makeatother
      \STATE\label{alg:gd:G:input}%
      \algorithmicrequire\ %
      \(\nabla\phi\colon G\to\frakg\), \(g_0\in G\); %
      \(\eta\colon\bbN_0\to\bbR\)%
      \STATE\label{alg:gd:G:ini}%
      \(x_0 \assign 0\)%
      \FOR{\(k = 1\) \TO \(N-1\)}\label{alg:gd:G:loop}%
      \STATE\label{alg:gd:G:step:gd}%
      \(x_{k+1} \assign x_k - \eta_k\nabla\phi(g_k)\)%
      \STATE\label{alg:gd:G:solve}%
      \(g_{k+1} \assign g_k\tau(\xi_k)\) such that %
      \(\xi_k = \diff\rtr_{\xi_k}^t\left( \Ad_{g_k}^t\Delta x_k \right)\)
      \ENDFOR\label{alg:gd:G:end}%
      \STATE \algorithmicensure\ \(g_N\)
    \end{algorithmic}
  \end{minipage}
\end{algorithm}


\section{Derivation}
\label{sec:devel}
In \Cref{sec:devel:direct}, we derive our novel scheme for Lie groups \cref{eq:mm:G}, which was previously introduced in \cref{alg:mm:G}. Later, in \Cref{sec:devel:pontryagin}, we demonstrate that this derivation can be obtained as a particular case of the Hamilton-Pontryagin framework developed in \cite{Bou1}. However, first, we shall recall the variational nature of PHB and NAG in the Euclidean case. For an introduction to variational integrators, we refer the reader to \cite{marsden-west}, and to \Cref{sec:euler-lagrange} for the case of Lie groups.

Classical and accelerated momentum methods, \emph{e.g.} Polyak's Heavy Ball and Nesterov's Accelerated Gradient, are equivalent to the discrete Euler-Lagrange equations of a particular discrete Lagrangian system on a path-connected open subset $D$ in the flat space $\bbR^n$ (confer with \cite{CMM23}). For a \(\calC^1\)-function \(\phi\colon D\subset\bbR^n\to\bbR\), these equations are
\begin{equation}
  \label{eq:DEL:Rn}
  \Delta x_k = \mu_k\Delta\left[ x_{k-1} - \varepsilon\eta_{k-1}\nabla\phi(x_{k-1})\right] - \eta_k\nabla\phi(x_k) \,,
\end{equation}
where \(\Delta\) is the forward difference operator, \(\mu_k\) and \(\eta_k\) are suitable coefficients (the method's strategy), and \(\varepsilon\) is a Boolean coefficient: \(\varepsilon=0\) for PHB and \(\varepsilon=1\) for NAG. The terms accompanying \(\varepsilon\) are associated to a force (as we will see later in the generalized framework of Lie groups), hence NAG is in fact PHB with forces.

This equation may be split in two steps to determine \(x_{k+1}\) from \(x_k\) and \(x_{k-1}\): a \emph{gradient (descent)} step \cref{eq:mm:Rn:gradient}, and a \emph{momentum} step \cref{eq:mm:Rn:momentum}:
\begin{subequations}
  \label{eq:mm:Rn}
  \begin{align}
    \label{eq:mm:Rn:gradient}
    y_{k+1} ={}& x_k - \eta_k\nabla \phi(x_k) \,,\\
    \label{eq:mm:Rn:momentum}
    x_{k+1} ={}& y_{k+1} + \mu_k\Delta z_k \,,
  \end{align}
\end{subequations}
where the variable \(z\) has a different meaning depending on the family of choice, \(z\equiv x_{.-1}\) for PHB, and \(z\equiv y\) for NAG. Equation \cref{eq:mm:Rn:gradient} should be viewed as an auxiliary definition that transforms \cref{eq:DEL:Rn} into \cref{eq:mm:Rn:momentum} and vice versa. Hence, although \(x\)'s and \(y\)'s follow a trajectory towards the argument minimum of \(\phi\), strictly speaking \(x_k\) is the natural one.

\subsection{A direct approach on Lie groups}
\label{sec:devel:direct}
We now derive \cref{alg:mm:G}, a novel class of methods on Lie groups, analogous to the classical PHB and NAG schemes. To this end, consider a real-valued \(\calC^1\) function \(\phi\) defined on a path-connected open subset \(D\) of a Lie group \(G\), that is, \(\phi\colon D \subseteq G \to \bbR\). Assume furthermore that \(\phi\) has a single local minimum in \(D\),
\begin{equation}
  \label{eq:optim}
  \textstyle
  \gopt = \argmin_{g\in D} \phi(g)\,.
\end{equation}

\begin{subequations}
  \label{eq:dsys}
  We define on \(D\times D\subset G\times G\) the discrete time-dependent Lagrangian system with forces \cite{CMM23,marsden-west}
  \begin{align}
\label{eq:dsys:lagrangian}
    l_k(w_0,w_1) \coloneq{} & a_k\tfrac12\|\rct(\Delta w_0)\|^2-b^-_k\phi(w_0)-b^+_{k+1}\phi(w_1) \,,\\
    \label{eq:dsys:forcem}
    f^-_k(w_0,w_1) \coloneq{} & {-\tfrac{a_{k-1}}{a_k}}(b^-_k+b^+_k)\diff\phi(w_0) \,,\\
    \label{eq:dsys:forcep}
    f^+_k(w_0,w_1) \coloneq{} & \phantom{-\tfrac{a_{k-1}}{a_k}}(b^-_k+b^+_k)\diff\phi(w_0)\circ\ra_{(\Delta w_0)^{-1}} \,,
  \end{align}
\end{subequations}
where \(a_k>0,b^\pm_k\) are arbitrary (but fixed) sequences of coefficients, $(w_0, w_1)\in D\times D$ and $\tau$ is a given retraction map (\Cref{sec:lie.groups}).
The discrete Euler-Lagrange equations of a free/forced system are (\Cref{sec:euler-lagrange}):
\[
  D_1l_{k+1}(w_1,w_2) + D_2l_k(w_0,w_1)
  + \varepsilon f^-_{k+1}(w_1,w_2) + \varepsilon f^+_k(w_0,w_1)
  = 0 \in\Tan_{w_1}^*G \,,
\]
where, as earlier, \(\varepsilon\) is a Boolean coefficient: \(\varepsilon=0\), free system; \(\varepsilon=1\), forced system.
Taking into account that
\[
  \pp[\rct(\Delta w_0)]{w_0} = -\Tan_{\Delta w_0}\rct\circ\la_{w_0^{-1}}\circ\ra_{\Delta w_0}
  \qquand
  \pp[\rct(\Delta w_0)]{w_1} = \Tan_{\Delta w_0}\rct\circ\la_{w_0^{-1}}\,,
\]
we obtain in this particular case
\begin{multline*}
 \begin{aligned}
  -a_{k+1}\ra_{\Delta w_1}^*\la_{w_1^{-1}}^*(\Tan_{\Delta w_1}\rct)^*&((\rct(\Delta w_1))^\flat)-b^-_{k+1}\diff\phi(w_1)\\
  {}+a_k\la_{w_0^{-1}}^*(\Tan_{\Delta w_0}\rct)^*&((\rct(\Delta w_0))^\flat)-b^+_{k+1}\diff\phi(w_1)
 \end{aligned}\\
  {}-\varepsilon\tfrac{a_k}{a_{k+1}}(b^-_{k+1}+b^+_{k+1})\diff\phi(w_1)
  +\varepsilon(b^-_k+b^+_k)\ra_{(\Delta w_0)^{-1}}^*\diff\phi(w_0)
  = 0\in\Tan_{w_1}^*G \,,
\end{multline*}
where \((\cdot)^\flat\) is the musical flat operator.
Divide by \(-a_{k+1}\), reorder terms, pull back to the identity by the right action, and apply the musical sharp operator \((\cdot)^\sharp\) to get
\begin{subequations}
  \label{eq:DEL:G:glob}
  \begin{equation}
    \label{eq:DEL:G}
    \Delta x_{k+1} = \mu_{k+1}\left(
      \Delta x_k
      - \varepsilon\Delta\left[\eta_k\nabla\phi(w_0)\right]
    \right)
    - \eta_{k+1}\nabla\phi(w_1)
    \in\frakg\,,
  \end{equation}
  where
  \[
    \mu_k\coloneq\tfrac{a_{k-1}}{a_k}\,,\quad
    \eta_k\coloneq\tfrac{b^-_k+b^+_k}{a_k}\,, \quand
    \Delta x_k \coloneq \left( \ra_{w_1}^*\la_{w_0^{-1}}^*(\Tan_{\Delta w_0}\rct)^*(\rct(\Delta w_0))^\flat \right)^\sharp\,.
  \]
  This last equation can be rewritten as
  \begin{equation}
    \label{eq:DEL:G:reconstruction}
    \Delta x_k = \left( \diff\rtr_{\rct(\Delta w_0)}^{-1}\circ\Ad_{w_0^{-1}}\right )^t\rct(\Delta w_0)\,.
  \end{equation}
  Indeed,
  \begin{align*}
    \Delta x_k ={}& \left( (\Tan_{\Delta w_0}\rct\circ\la_{w_0^{-1}}\circ\ra_{w_1})^*(\rct(\Delta w_0))^\flat \right)^\sharp\\
    ={}& \left( \Tan_{\Delta w_0}\rct\circ\la_{w_0^{-1}}\circ\ra_{w_1} \right)^t\rct(\Delta w_0)\\
    ={}& \left( \diff\rtr_{\rct(\Delta w_0)}^{-1}\circ\ra_{(\Delta w_0)^{-1}}\circ\la_{w_0^{-1}}\circ\ra_{w_1} \right)^t\rct(\Delta w_0)\,,
  \end{align*}
  where we have first used a simple relation between the musical operators, the dual map, and the map transpose, \((A^*v^\flat)^\sharp=A^tv\), then the definition of \(\tau\)'s right-trivialized tangent \cref{eq:retraction:tangent}, and finally the commutativity of the left and right actions to get the adjoint representation after simplification.
\end{subequations}

The set of equations in \cref{eq:DEL:G:glob} defines two families of methods---or, equivalently, a family of twin methods---which we refer to as momentum methods for Lie groups: the classical variant when \(\varepsilon = 0\), and the accelerated variant when \(\varepsilon = 1\). Although \cref{eq:DEL:G} is formally identical to its Euclidean counterpart \cref{eq:DEL:Rn}, for the time being, it cannot be expressed in the form of \cref{eq:mm:Rn}. In \cref{eq:DEL:G:glob}, solely the bracketing \(\Delta\) corresponds to the usual difference operator, while \(\Delta w_k\) represents the group right-transition, and \(\Delta x_k\) is merely suggestive notation. That is, there is no canonical choice of \(x_k\) and \(x_{k+1}\) such that \(\Delta x_k = x_{k+1} - x_k\), which prevents the introduction of \cref{eq:mm:Rn:gradient} to rewrite \cref{eq:DEL:G} in the form of \cref{eq:mm:Rn:momentum}. However, if we set \(x_0\) to any fixed value (for instance, \(x_0=0\in\frakg\)), then all \(x_{k+1}=x_k+\Delta x_k\) become defined recursively.

We may now rewrite \cref{eq:DEL:G:glob} for \(w_j=g_{k+j}\) in the form of \cref{eq:mm:Rn}:
\begin{subequations}
  \label{eq:mm:G}
  \begin{align}
    \label{eq:mm:G:gradient}
    y_{k+1} ={}& x_k - \eta_k\nabla\phi(g_k) \,,\\
    \label{eq:mm:G:z}
    z_{k+1} ={}& (1-\varepsilon)x_k + \varepsilon y_{k+1}\,,\\
    \label{eq:mm:G:momentum}
    x_{k+1} ={}& y_{k+1} + \mu_k\Delta z_k \,,\\
    \label{eq:mm:G:reconstruction}
    g_{k+1} ={}& g_k\Delta g_k \quad\text{such that}\quad  \rct(\Delta g_k) = \diff\rtr_{\rct(\Delta g_k)}^t\left( \Ad_{g_k}^t\Delta x_k \right) \,,
  \end{align}
\end{subequations}
where \cref{eq:mm:G:z} has been added for convenience, and where \cref{eq:mm:G:reconstruction} is the \emph{reconstruction} step from Equation \cref{eq:DEL:G:reconstruction}. Note that this equation is implicit. In fact, \(\xi_k\coloneqq\rct(\Delta g_k)\) is a solution of the fixed point equation
\(
  \xi = \diff\rtr_\xi^t\eta
\)
with \(\eta\coloneqq\Ad_{w_0}^t\Delta x_k\).

As far as we know, Eqs.~\cref{eq:DEL:G,eq:mm:G} constitute novel formulations of classical and accelerated momentum methods on Lie groups. The computational cost is primarily concentrated in the gradient evaluation in \cref{eq:mm:G:gradient}, and partially in the reconstruction of group elements via \cref{eq:mm:G:reconstruction}. However, in certain cases (\cref{sec:solvers}), this equation turns out to be explicit, lowering the computational burden.

Being \cref{eq:DEL:G} a difference equation of order 2, two initial values \(g_0,g_1\in G\) sufficiently close to \(\gopt\) are required. Given \(g_0\), take \(g_1=g_0\), for which \cref{eq:DEL:G:reconstruction} gives \(\Delta x_0=0\in\frakg\). Then define \(y_1\) and \(z_1\) using Equations \cref{eq:mm:G:gradient,eq:mm:G:z} with \(k=0\), before running the whole scheme \cref{eq:mm:G} for \(k\geq1\). In summary,
\begin{equation}
  \label{eq:mm:G:inival}
  g_1 = g_0 \,,\ \
  \Delta x_0 = 0 \,,\ \
  (x_0 = 0) \,,\ \
  y_1 = x_0 - \eta_0\nabla\phi(g_0) \,,\ \
  z_1 = (1-\varepsilon)x_0 + \varepsilon y_0 \,.
\end{equation}

On a side note, there is a workaround to avoid having to set \(x_0\): Subtract two consecutive sets of Eqs. \cref{eq:mm:G} to get
\begin{subequations}
  \label{eq:mm:G:double}
  \begin{align}
    \label{eq:mm:G:double:gradient}
    \Delta y_{k+1} ={}& \Delta x_k - \Delta[\eta_k\nabla\phi(g_k)] \,,\\
    \label{eq:mm:G:double:z}
    \Delta z_{k+1} ={}& (1-\varepsilon)\Delta x_k + \varepsilon \Delta y_{k+1}\,,\\
    \label{eq:mm:G:double:momentum}
    \Delta x_{k+1} ={}& \Delta y_{k+1} + \Delta[\mu_k\Delta z_k] \,,\\
    \label{eq:mm:G:double:reconstruction}
    g_{k+1} ={}& g_k\Delta g_k \quad\text{such that}\quad  \rct(\Delta g_k) = \diff\rtr_{\rct(\Delta g_k)}^t\left( \Ad_{g_k}^t\Delta x_k \right) \,.
  \end{align}
\end{subequations}
Although this does not increase significantly the overall cost, its implementation would be slightly more cumbersome.

It is worth noting that in Eqs.~\cref{eq:mm:G,eq:mm:G:double}, the trivialization has not been explicitly stated. The same choice, whether right or left trivialization, must be made in Eqs.~\cref{eq:mm:G:gradient,eq:mm:G:reconstruction}, or in their doubled version, Eqs.~\cref{eq:mm:G:double:gradient,eq:mm:G:double:reconstruction}.

A final remark on the choice of strategy \((\mu_k, \eta_k)\). As in the Euclidean case (see \cite{CMM23}), and shown in the above discussion, these coefficients are linked to the Lagrangian parameters \((a_k, b_k)\) in \cref{eq:dsys}. Different choices lead to different convergence rates (cf.~\cite{SuBoCa16,WiWiJo16}), a topic that lies beyond the scope of the present work. It is worth noting, however, that since the Lagrangian parameters must be strictly or explicitly time-dependent, the strategy coefficients \((\mu_k, \eta_k)\) should, in principle, also vary with time. Nevertheless, there exist choices of Lagrangian coefficients for which the resulting strategy becomes constant (\Cref{sec:experiments}).

\subsection{The Hamilton-Pontryagin approach for Lie groups}
\label{sec:devel:pontryagin}
{Mo\-men\-tum-\allowbreak ba\-sed} gradient descent methods on Lie groups, such as Eq.~\cref{eq:mm:G}, can also be derived from a Hamilton-Pontryagin variational principle, yielding dynamical equations similar to those in \cite{Bou1}. Given the relevance of this approach, we devote this section to the derivation of forward and backward explicit Euler methods on Lie groups. We also show that the previous derivation can be interpreted as a particular instance of this general framework.

Let \(\bar{l}\colon\bbZ\times G\times\frakg\to\bbR\) be a discrete time-dependent trivialized Lagrangian, define the discrete Lagrangian in Hamilton-Pontryagin form
\[
  \tilde{l}_k(z_k,z_{k+1}) \coloneqq \bar{l}_k(g_k,\xi_k) + \langle p_k, \rct(\Delta g_k)-\xi_k \rangle \,,
\]
where \(z_k=(g_k,\xi_k,p_k)\in G\times\frakg\times\frakg^*\).
The DEL equations \cref{eq:DEL:point-point} for such a Lagrangian read
\[
  \left\langle D_1\tilde{l}_k(z_k,z_{k+1}) + D_2\tilde{l}_{k-1}(z_{k-1},z_k) ,
    \delta z_k \right\rangle = 0
\]
for any variation \(\delta z_k\).
This equation decomposes with respect to \(\delta z_k=(\delta g_k,\delta\xi_k,\delta p_k)\) into
\begin{align*}
  \delta g\colon\ %
  &p_k\!\circ\!\Tan_{\Delta g_k}\rct\!\!\circ\!\la_{g_k^{-1}}
    +\partial_g\bar{l}_{k+1}(g_{k+1},\xi_{k+1})
    -p_{k+1}\!\circ\!\Tan_{\Delta g_{k+1}}\rct\!\!\circ\!\la_{g_{k+1}^{-1}}\!\circ\!\ra_{\Delta g_{k+1}}=0,\\
  \delta \xi\colon\ &\partial_\xi\bar{l}_k(g_k,\xi_k)-p_k=0,\\
  \delta p\colon\ &\rct(\Delta g_k)-\xi_k=0,
\end{align*}
which, after pulling the first to the identity, time-shifting the second, and using the definitions of \Cref{sec:lie.groups}, is rewritten in the form
\begin{subequations}
  \label{eq:DEL:G:pontryagin}
  \begin{align}
    \label{eq:DEL:G:pontryagin:reconstruction}
    g_{k+1} ={}& g_k \rtr(\xi_k) \,,\\
  \label{eq:DEL:G:pontryagin:momenta}
    p_{k+1} ={}& \partial_\xi\bar{l}_{k+1}(g_{k+1},\xi_{k+1}) \,,\\
  \label{eq:DEL:G:pontryagin:dynamics}
    \big(\diff\rtr_{\xi _{k+1}}^{-1}\big)^*p_{k+1}
    ={}& \Ad_{\Delta g_k}^*\big(\diff\rtr_{\xi _k}^{-1}\big)^*p_k
         + \la_{g_{k+1}}^*\partial_g\bar{l}_{k+1}(g_{k+1},\xi_{k+1}) \,.
  \end{align}
\end{subequations}

Had we defined the discrete Lagrangian in this alternate form
\[
  \tilde{l}_k(z_k,z_{k+1}) \coloneqq \bar{l}_k(g_k,\xi_k) + \langle p_{k+1}, \rct(\Delta g_k)-\xi_{k+1} \rangle \,,
\]
we would have ended up with the variational integrator
\begin{subequations}
  \label{eq:DEL:G:pontryagin:backward}
  \begin{align}
    \label{eq:DEL:G:pontryagin:backward:reconstruction}
    g_{k+1} ={}& g_k \rtr(\xi_{k+1}) \,,\\
    \label{eq:DEL:G:pontryagin:backward:momenta}
    p_{k+1} ={}& \partial_\xi\bar{l}_{k+1}(g_{k+1},\xi_{k+1}) \,,\\
    \label{eq:DEL:G:pontryagin:backward:dynamics}
    \big(\diff\rtr_{\xi_{k+1}}^{-1}\big)^*p_{k+1}
    ={}& \Ad_{\Delta g_k}^*\big(\diff\rtr_{\xi_k}^{-1}\big)^*p_k
         + \la_{g_k}^*\partial_g\bar{l}_k(g_k,\xi_k) \,.
  \end{align}
\end{subequations}

Assume the Lagrangian is left-invariant, that is, \(\partial_g\bar{l}_k=0\in\frakg^*\), and redefine the momenta as \(P_k\coloneq\big(\diff\rtr_{\xi_k}^{-1}\big)^*p_k\), then the scheme \cref{eq:DEL:G:pontryagin} can be rewritten as follows
\begin{subequations}
  \begin{align*}
    g_{k+1} ={}& g_k \rtr(\xi_{k+1}) \,,\\
    P_{k+1} ={}& \Ad_{\Delta g_k}^*P_k \,,\\
    \label{eq:DEL:G:pontryagin:momenta:invariant}
    P_{k+1} ={}& \big(\diff\rtr_{\xi_{k+1}}^{-1}\big)^*\partial_\xi\bar{l}_{k+1}(g_{k+1},\xi_{k+1}) \,,
  \end{align*}
\end{subequations}
which is explicit except for the last equation. Similarly, the scheme \cref{eq:DEL:G:pontryagin:backward} is backward (in time) explicit.
In fact, Equations \cref{eq:DEL:G:pontryagin,eq:DEL:G:pontryagin:backward} correspond, respectively, to Euler forward and backward methods (compare with \cite{Bou1}, Eqs.~(4.19) and (4.20), which are similar but slightly different).

For the particular case of the Lagrangian \cref{eq:dsys:lagrangian} or, rather, for the trivialized Lagrangian
\[ \bar{l}_k(g,\xi) \coloneqq l_k(g,g\mkern1mu\rtr(\xi)) = a_k\tfrac12\|\xi\|^2-b^-_k\phi(g)-b^+_{k+1}\phi(g\mkern1mu\rtr(\xi)) \,, \]
equation \cref{eq:DEL:G:pontryagin:dynamics} is equivalent to Equation \cref{eq:DEL:G} (with \(\varepsilon=0\)). To see this, simply compute the differential maps
\begin{align*}
  \partial_g\bar{l}_k(g_0,\xi_0)
  ={}& -b_k^-\diff\phi(g_0)-b_{k+1}^+\diff\phi(g_1)\circ\ra_{\Delta g_0}\\
  ={}& -b_k^-\ra_{g_0^{-1}}^*\nabla\phi(g_0)^\flat-b_{k+1}^+\ra_{g_0^{-1}}^*\nabla\phi(g_1)^\flat \,,\\
  \partial_\xi\bar{l}_k(g_0,\xi_0)
  ={}& a_k\xi_0^\flat-b_{k+1}^+\diff\phi(g_1)\circ\la_{g_0}\circ \Tan_{\xi_0}\rtr\\
  ={}& a_k\xi_0^\flat-b_{k+1}^+\big(\diff\rtr_{\xi_0}\big)^*\Ad_{g_0}^*\nabla\phi(g_1)^\flat \,,
\end{align*}
written in terms of the gradient and the trivialized tangent, and take into account the ``definitions'' \cref{eq:DEL:G:pontryagin:momenta,eq:DEL:G:reconstruction}.


\section{Examples}
\label{sec:examples}
For the numerical experiments presented in the next section, we consider combinations of different solutions to the reconstruction equation \cref{eq:mm:G:reconstruction} and various objective functions introduced herein. The examples are defined on the group of spatial rotations \(SO(3)\), consisting of orthogonal matrices with positive determinant. Its Lie algebra, \(\mathfrak{so}(3)\), is the space of skew-symmetric matrices. Throughout this section, \(R\) denotes a rotation matrix, while \(\Delta x\), \(\hat{\Omega}\), and \(\hat{\Theta}\) represent skew-symmetric matrices, with the latter serving a utilitarian role in the discussions below. Additionally, \((\cdot)^-\) denotes the skew-symmetric part of a matrix, and \((\cdot)^\wedge\) denotes the representation of a vector in \(\mathbb{R}^3\) as a skew-symmetric matrix.

\subsection{Solvers for the reconstruction equation\texorpdfstring{\ \cref{eq:mm:G:reconstruction}}{}}
\label{sec:solvers}
Natural or common retraction maps on \(\frakso(3)\) are the matrix exponential (\Cref{sec:expm}) and the Cayley transform (\Cref{sec:cayley}). Another case of interest is the skewsymmetric part of a rotation \(R\in SO(3)\), the inverse of a certain retraction (\Cref{sec:skew}).

Before we define specific objective functions to be optimized, we must first observe that Equation \cref{eq:mm:G:reconstruction} can be rendered explicit for the selected retractions. In fact, in these cases, it is equivalent to the following expressions
\begin{subequations}
  \label{eq:solvers}
  \begin{align}
    \label{eq:solver:exp}
    R_{k+1} ={}& R_k\exp(\Delta x_k) \,,&
    R_{k+1} ={}& \exp(\Delta x_k)R_k \,,\\
    \label{eq:solver:cay}
    R_{k+1} ={}& R_k\cay(2\lambda\Delta x_k) \,,&
    R_{k+1} ={}& \cay(2\lambda\Delta x_k)R_k \,,\\
    \label{eq:solver:skew}
    R_{k+1} ={}& R_k\unskew(\gamma\Delta x_k) \,,&
    R_{k+1} ={}& \unskew(\gamma\Delta x_k)R_k \,,
  \end{align}
\end{subequations}
where \(R_k \in SO(3)\) and \(\Delta x_k\in\frakso(3)\equiv\bbR^3\), and where the side in which the equations appear has a direct correspondence with the choice of left or right acting group transitions. Besides, for the Cayley transform the coefficient \(\lambda\) is given by
\[
  \lambda = \frac1{1+(\Lambda-\frac1{3\Lambda})^2}
  \quad\text{with}\quad
  \Lambda = \sqrt[3]{\|\Delta x_k\|+\sqrt{\|\Delta x_k\|^2+\tfrac1{27}}}
  \,,
\]
and for the inverse skewsymmetric projection the coefficient \(\gamma\) is a solution to
\[ \|\Delta x_k\|^2\gamma^4-2\gamma+1 = 0 \,, \]
whose solution is unique for \(\|\Delta x_k\|<1\).

Indeed, Equation \cref{eq:DEL:G:reconstruction} is equivalent to
\[
  \left( \diff\rtr_{\rct(\Delta g_k)}^{-1} \right)^t\left( \rct(\Delta g_k) \right) = \Ad_{g_k}^t\left( \Delta x_k \right) \,.
\]
For \(\tau=\exp\colon SO(3)\to\frakso(3)\), this relation reads
\[
  \dlog(\log(\Delta R_k))^t\left(\log(\Delta R_k)\right) = \Ad_{R_k}^t(\Delta x_k) \,.
\]
Since \(\dlog(\hat x)^t(x)=x\) (confer with Equation \cref{eq:exp:dinv}), we get
\[
  \log(\Delta R_k) = R_k^{\,t}\Delta x_kR_k \,,
\]
which finally gives Equation \cref{eq:solver:exp} (right, left is analogous),
\[
  R_{k+1} = R_k\exp(R_k^{\,t}\Delta x_kR_k) = \exp(\Delta x_k)R_k
\]
where we have used the fact that the exponential map commutes with conjugation.

The case for the Cayley transform is slightly different, since now \(\dcay^{-1}(\hat x)^t(x)=\tfrac{1+\|x\|^2}2x\)  (confer with Equation \cref{eq:cay:dinv}). We have still
\[
  \tfrac{1+\|\Omega_k\|^2}2\Omega_k
  = \dcay^{-1}(\hat\Omega_k)^t(\Omega_k)
  = R_k^{\,t}\Delta x_kR_k \,,
\]
where \(\hat\Omega_k=\cay^{-1}(\Delta R_k)\).
Applying the norm to both sides results in
\[
  \|\Omega_k\|^3+\|\Omega_k\|-2\|\Delta x_k\| = 0 \,,
\]
a third order algebraic equation for \(\|\Omega_k\|\) with a single real root,
\(\Lambda-\tfrac1{3\Lambda}\),
as in Equation \cref{eq:solver:cay}, which is proven once the commutation between the Cayley transform and the conjugation is taken into account,
\[
  R_{k+1}
  = R_k\cay\left( \tfrac2{1+\|\Omega_k\|^2} R_k^{\,t}\Delta x_kR_k \right)
  = \cay\left( \tfrac2{1+\|\Omega_k\|^2}\Delta x_k\right)R_k \,.
\]

Finally, an analogous derivation follows for the \(\unskew\) retraction map. In fact, using the relation \(\dskew(\hat x)^t(x) = \gamma^{-1}x\) (confer with Equations \cref{eq:skew:d}), expression \cref{eq:DEL:G:reconstruction} gives in this case
\[
  \gamma^{-1}\Omega_k
  = \dskew(\hat\Omega_k)^t(\Omega_k)
  = R_k^{\,t}\Delta x_kR_k \,,
\]
where \(\hat\Omega_k=\skew(\Delta R_k)\), and \(\gamma^{-1}=1+\sqrt{1-\|\Omega_k\|^2}\).
As in the previous cases, commutation gives the result.
However, prior to that, \(\gamma\) should be determined. Taking norms on both sides gives \(\gamma^{-1}\|\Omega_k\| = \|\Delta x_k\|\), which means \(\gamma\) is a solution of
\[
  2\gamma^{-3}-\gamma^{-4}=\|\Delta x_k\|^2 \,.
\]
Explicit solutions to this equation for \(\gamma^{-1}\) in function of \(\|\Delta x_k\|\) may be given, however these solutions have a rather involved expression. An alternative is to use a nonlinear solver such as Newton-Raphson starting from a safe initial guess. From \(\gamma\)'s definition \(\tfrac12\leq\gamma\leq1\), and the solution to the above equation is unique for \(\tfrac12\leq\gamma\leq\tfrac23\). So a safe initial guess is \(\gamma_0=\tfrac7{12}\approx0.583\).

\subsection{Objective functions}
\label{sec:objectives}
We consider four different functions defined by restriction or retraction.

\subsubsection{Restricted squared Frobenius norm}
Consider the Frobenius (or entrywise) norm on the space of squared matrices \(\mathcal{M}_{3\times3}(\bbR)\) and let \(f\colon\mathcal{M}_{3\times3}(\bbR)\to\bbR\), \(A\mapsto\tfrac12\norm{A-I}^2\). We define
\begin{equation}
  \label{eq:objective:frobenius}
  \phi\coloneq f|_{SO(3)} \text{, for which } \nabla\phi(R) = R^- \,.
\end{equation}
Since \(f\) is continuous and \(SO(3)\subset\mathcal{M}_{3\times3}(\bbR)\) is compact, we know that \(\phi\) attains its global minimum and maximum values (\(0\) and \(4\)), which in fact occurs respectively at the identity \(I\) and at rotations with \(-1\) trace.

\subsubsection{Restricted Rosenbrock function}
Rosenbrock's function \cite{Ro60}, whose expression is
\[ \ros(x,y;a,b) = (a-x)^2 + b(y-x^2)^2 \,, \] with parameters \(a,b>0\), represents a banana-shaped flat-valley surrounded by steep walls with a unique critical point and global minimum at \(a,a^2\), whose search by numerical means is difficult, hence its use to test and benchmark optimizers. We consider here its generalization to higher dimensions, \(n>2\), namely
\begin{equation}
  \label{eq:rosenbrock}
  \ros(x) = \sum_{i=1}^{n-1}\ros(x_i,x_{i+1};1,100) = \sum_{i=1}^{n-1}\left[ (1-x_i)^2+100(x_{i+1}-x_i^2)^2 \right] \,.
\end{equation}
As in the two-dimensional case, the function has a global minimum at \((1,1,\dots,1)\) but, unlike it, also has a local minimum close to \((-1,1,\dots,1)\) (the higher is the dimension, the closer it gets).

Consider the function \(g\colon\mathcal{M}_{3\times3}(\bbR)\to\bbR\), \(A\mapsto\ros(\mathbf{1}+A-I)\), where \(\mathbf{1}\) is a matrix filled with \(1\)'s, and where the entries of the matrix to apply the Rosenbrock function ought to be taken columnwise. We define by restriction
\begin{equation}
  \label{eq:objective:rosenbrock:restricted}
  \phi\coloneq g|_{SO(3)} \text{, for which } \nabla\phi(R) = (R\cdot\nabla\ros(\mathbf{1}+R-I))^- \,.
\end{equation}
The unique global minimum is attained at the identity \(I\) and is surrounded by other local minima. The global maximum is presumably at \(\tfrac13(2\mathbf{1}-3I)\).

\subsubsection{Retracted Rosenbrock function}
\begin{subequations}
\label{eq:objective:rosenbrock}
As objective function, we consider a composition of either of the chosen retractions with the Rosenbrock function in \(\bbR^3\), that is,
\begin{equation}
  \label{eq:objective:rosenbrock:fun}
  \phi(R) \coloneq \ros(\rct(R)^\vee) \,,\ \forall R\in SO(3) \,,
\end{equation}
where \(\ros(x,y,z)=(1-x)^2+100\cdot(y-x^2)^2+(1-y)^2+100\cdot(z-y^2)^2\). It is then readily seen that the objective function \(\phi\) has a unique global minimum at \(\tau\big(\widehat{(1,1,1)}\big)\).

To compute \(\nabla\phi\), given \(R\in SO(3)\), let \(\hat\Omega=\rct(R)\in\frakso(3)\), and take \(\hat\Theta\in\frakso(3)\) arbitrary, then we get
\begin{align*}
  \left\langle \nabla\phi(R), \hat\Theta \right\rangle
  ={}& \diff\phi(R)(\hat\Theta R)\\
  ={}& \dros(\Omega)\cdot(\Tan_{\rtr(\hat\Omega)}\rct(\hat\Theta R))^\vee\\
  ={}& \dros(\Omega)\cdot((\Tan_{\hat\Omega}\rtr)^{-1}(\hat\Theta R))^\vee\\
  ={}& \dros(\Omega)\cdot\drct(\hat\Omega)(\Theta)\\
  ={}& \left\langle \nabla\!\ros(\Omega), \drct(\hat\Omega)(\Theta) \right\rangle\\
  ={}& \left\langle \big(\drct(\hat\Omega)\big)^t\cdot\nabla\!\ros(\Omega), \Theta \right\rangle \,.
\end{align*}
where we have applied (in order) the trivialized gradient definition, the chain rule, the inverse function theorem, the trivialized tangent definition, the regular gradient definition, and the linear map transposition.
Therefore, for the particular cases of the exponential and the Cayley transform (see \cref{eq:exp:dinv,eq:cay:dinv}), we have
\begin{equation}
  \label{eq:objective:rosenbrock:grad}
  \nabla\phi(R) = \big(\drct(-\hat\Omega)\cdot\nabla\!\ros(\Omega)\big)^\wedge \,.
\end{equation}
\end{subequations}


\section{Experiments and results}
\label{sec:experiments}
Several experiments have been conducted, and we present but a meaningful subset in the following figures. The plots illustrate on a logarithmic scale the residue of the objectives functions described in the preceding \Cref{sec:objectives}.
We consider three optimization methods: gradient descent (orange), Polyak's heavy ball (blue), and Nesterov's accelerated gradient (green).
For reference, we include sequences of the form \(O(1/k^2)\) (red). Our exploration involves the three solvers of \Cref{sec:solvers}: one (left) is based on the exponential map, Eq.~\cref{eq:solver:exp}; another (mid) employs the Cayley transform, Eq.~\cref{eq:solver:cay}; the third (right) uses the inverse of the skewsymmetric projection,  Eq.~\cref{eq:solver:skew}.
The chosen strategies \((\mu_k, \eta_k)\) vary across experiments but are constant in each case, and are (approximately) derived from an exponentially dilated Lagrangian (cf.~\cite{CMM23}). A strategy is considered more \emph{aggressive} (resp., \emph{conservative}) than another if one or both of its coefficients are larger (resp., smaller).

The experiments were implemented in Julia \cite{bezanson2017julia,julia1070} and are available at an open access repository \cite{cmcxyz}. They only pretend to show that, in general, the schemes perform as expected, but do not for particular cases that we highlight. This shows that a numerical analysis, out of the scope of the present paper, is nonetheless of interest.

\newcommand\mycaption[5]{\caption{Residue (log scale) \emph{vs.} epoch for #1. Simulation run for #2 epochs from initial guess \(R_0=#3\) with constant strategy \(\mu_0=#4\) (0 for GD), \(\eta_0=#5\).}}
\begin{figure}
  \centering
  \includesvg[width=.325\textwidth]{figures/frobenius1-exp}
  \includesvg[width=.325\textwidth]{figures/frobenius1-cay}
  \includesvg[width=.325\textwidth]{figures/frobenius1-skw}
  \mycaption{the restricted Frobenius norm}{100}{\cay(1,1,1)}{0.7}{0.1}
  \label{fig:frobenius1}
  \bigskip
  \includesvg[width=.325\textwidth]{figures/frobenius2-exp}
  \includesvg[width=.325\textwidth]{figures/frobenius2-cay}
  \includesvg[width=.325\textwidth]{figures/frobenius2-skw}
  \mycaption{the restricted Frobenius norm}{250}{\cay(1,1,1)}{0.7}{0.01}
  \label{fig:frobenius2}
\end{figure}
\begin{figure}
  \centering
  \includesvg[width=.325\textwidth]{figures/rosenbrock91-exp}
  \includesvg[width=.325\textwidth]{figures/rosenbrock91-cay}
  \includesvg[width=.325\textwidth]{figures/rosenbrock91-skw}
  \mycaption{the restricted Rosenbrock function}{100}{\cay(0.1,0.1,0.1)}{0.25}{0.0001}
  \label{fig:rosenbrock91}
  \bigskip
  \centering
  \includesvg[width=.325\textwidth]{figures/rosenbrock92-exp}
  \includesvg[width=.325\textwidth]{figures/rosenbrock92-cay}
  \includesvg[width=.325\textwidth]{figures/rosenbrock92-skw}
  \mycaption{the restricted Rosenbrock function}{100}{\cay(0.1,0.1,0.1)}{0.7}{0.0001}
  \label{fig:rosenbrock92}
\end{figure}
\begin{figure}
  \centering
  \includesvg[width=.325\textwidth]{figures/rosenbrock3exp-exp}
  \includesvg[width=.325\textwidth]{figures/rosenbrock3exp-cay}
  \includesvg[width=.325\textwidth]{figures/rosenbrock3exp-skw}
  \mycaption{the Rosenbrock function retracted by \(\tau=\exp\)}{1000}{\exp(0,0,1)}{0.99}{0.0001}
  \label{fig:rosenbrock31}
  \bigskip
  \includesvg[width=.325\textwidth]{figures/rosenbrock3cay-exp}
  \includesvg[width=.325\textwidth]{figures/rosenbrock3cay-cay}
  \includesvg[width=.325\textwidth]{figures/rosenbrock3cay-skw}
  \mycaption{the Rosenbrock function retracted by \(\tau=\cay\)}{1000}{\cay(0,0,1)}{0.99}{0.0001}
  \label{fig:rosenbrock32}
\end{figure}

As expected, we observe that in all cases the three methods outperform the reference rate \(O(1/k^2)\). In most experiments, NAG achieves the best performance, followed by PHB, with GD being the least effective. However, this hierarchy does not always hold and depends on the chosen solver or strategy. We have deliberately retained such cases in the presentation to allow for further discussion.

For instance, in \cref{fig:frobenius1} (mid), when optimizing the squared Frobenius norm using a ``medium-large'' momentum coefficient combined with a ``large'' learning rate (i.e., a large time step) and the Cayley transform, GD outperforms clearly PHB and slightly NAG. However, this is no longer the case when a more conservative strategy is adopted in terms of the learning rate, as shown in \cref{fig:frobenius2}, where the expected hierarchy is recovered. Observe also that, under this latter strategy, PHB slightly outperforms NAG when using the exponential and skew-based solvers. Both figures correspond to the simplest objective function considered: the squared Frobenius norm.

In the case of the 9-dimensional Rosenbrock function, depicted in \cref{fig:rosenbrock91}, a conservative strategy with both momentum and learning rate set to small values yields results similar to those in \cref{fig:frobenius2}, where GD proves to be the least effective. In contrast, and more notably, a more aggressive strategy in terms of momentum reduces the performance of PHB and NAG compared to GD when using the Cayley-based solver, \Cref{fig:rosenbrock92} (mid), while the expected ranking is preserved for the other two solvers (left and right).

For both versions of the retracted Rosenbrock function, \Cref{fig:rosenbrock31,fig:rosenbrock32}, we observe the expected behavior: a clear improvement when using momentum-based methods over GD. Additionally, the momentum-based optimization trajectories exhibit a characteristic circling pattern around the minimizer, reminiscent of a ball rolling inside a bowl.


\section{Conclusions}
\label{sec:conclusions}
We present a variational derivation of first-order momentum methods for Lie groups. These schemes generalize the well-known PHB and NAG methods in \(\mathbb{R}^n\). These familiar methods emerge as special cases when considering the group of translations in \(\mathbb{R}^n\) with the identity as the retraction map. In fact, the methods applied to both Euclidean space and Lie groups share a common formal structure, Eqs.~\cref{eq:DEL:Rn,eq:DEL:G}, albeit with few distinctions. As in general, a Lie group is not a linear space, we cannot write the right translation element \(\Delta g_0=g_0^{-1}g_1\) as the difference \(g_1-g_0\). To address this, we resort to pull the problem to the Lie algebra associated with the group. The intricate relationship between the group and its algebra is captured by a novel equation, termed the reconstruction equation, Eq.~\cref{eq:DEL:G:reconstruction}. Apart from this equation, the schemes are explicit, and in specific scenarios, this equation can also be rendered explicit, Eqs.~\cref{eq:solvers}, thus reducing in principle the overall computational cost. Notably, this holds true for the exponential map, the Cayley transform, and the inverse of the skew-symmetric projection. Furthermore, our method can be implemented either directly in terms of \(x_k\) by setting \(x_0 = 0\), \cref{alg:mm:G}, or in terms of \(\Delta x_k\) using an overlapped approach, Eqs.~\cref{eq:mm:G:double}.

The methods have been formulated by exploiting the inherent geometrical structure of these spaces. They are equivalent to the Euler-Lagrange equations of specific Lagrangian systems, Eq.~\cref{eq:dsys}. In addition, these methods admit an alternative formulation in Hamilton-Pontryagin form, which connects them to the forward and backward Euler methods, Eqs.~\cref{eq:DEL:G:pontryagin,eq:DEL:G:pontryagin:backward}. This twofold derivation of \cref{alg:mm:G} provides a form of theoretical validation for the proposed scheme.

In general, numerical results align with expectations, \cref{fig:frobenius2,fig:rosenbrock91,fig:rosenbrock31,fig:rosenbrock32}. However, there exist cases that deviate from this general trend, \cref{fig:frobenius1,fig:rosenbrock92} (both mid), highlighting the need for a more detailed numerical analysis, which lies beyond the scope of the present work. Such an analysis should into account not only the properties of the objective function, the scheme's family, and the chosen strategy, but also the geometric aspects of the Lie group, as conveyed through the retraction map.

\appendix

\section{Retractions on Lie groups}
\label{sec:lie.groups}
Let \(G\) be a Lie group, \(TG\) denotes the tangent bundle, \(\frakg=\Tan_eG\) its Lie algebra, where \(e\) is the neutral element of \(G\), and \(T^*G\) and \(\frakg^*\) their duals. The left and right actions (or translations) of the group are denoted \(\la_g\) and \(\ra_h\), respectively, so that \(\la_g(h)=gh=\ra_h(g)\). It readly seen that left and right translation commute, that is, \(\la_g\circ\ra_h=\ra_h\circ\la_g\). Moreover, these maps allow for the trivialization of the tangent and cotangent bundles. For the left action:
\begin{align*}
  TG&\To G\times \frakg&
  T^*G&\To G\times \frakg^*\\
  (g,\dot g)&\longmapsto \big(g, \Tan_g\la_{g^{-1}}\dot g\big)&
  (g, \alpha)&\longmapsto \big(g, (\Tan_e\la_{g})^*\alpha\big)
\end{align*}
Analogously for the right action.

The \emph{conjugation} is the map \(\mathrm{C}_g\coloneqq\la_g\circ\ra_{g^{-1}}\colon G\to G\), the \emph{adjoint group representation} is \(\Ad\colon G\to Gl(\frakg)\) such that \(\Ad_g\coloneqq \Tan_e\mathrm{C}_g\colon\frakg\to\frakg\), and the \emph{adjoint algebra representation} is \(\ad\coloneqq \Tan_e\Ad\colon\frakg\to\frakgl(\frakg)\) so that \(\ad_\xi\eta=[\xi,\eta]\).

A \emph{retraction} on \(G\) is a mapping \(\tau\colon\frakg\to G\), which is an analytic local diffeomorphism around the identity such that \(\tau(\xi)\tau(-\xi)=e\) for any \(\xi\in\frakg\). Thereby, \(\tau\) provides a local chart on the Lie group. A particular case of retraction is the exponential map.

Given a retraction \(\tau\colon\frakg\to G\), we define its \emph{right-trivialized tangent} \cite{Bou1} as the mapping \(\diff\tau\colon\frakg\times\frakg\to \frakg\) given for any \(\xi\in\frakg\) by
\begin{equation}
  \label{eq:retraction:tangent}
  \diff\tau(\xi,\cdot) = \diff\tau_\xi \coloneqq \Tan_g\ra_{g^{-1}}\circ\:\Tan_\xi\tau \,,
\end{equation}
where \(g=\tau(\xi)\), therefore \(g^{-1} = \tau(\xi)^{-1} = \tau(-\xi)\).
The \emph{right-trivialized inverse tangent} of \(\tau\) is the mapping \(\diff\tau^{-1}\colon\frakg\times\frakg\to\frakg\)
\begin{equation}
  \label{eq:retraction:tangent:inverse}
  \diff\tau^{-1}(\xi,\cdot) = \diff\tau^{-1}_\xi \coloneqq (\diff\tau_\xi)^{-1} = \Tan_g\tau^{-1}\circ\:\Tan_e\ra_g\,.
\end{equation}
The left-trivialized direct and inverse tangent are defined analogously.

The trivialized tangents have a simple relation with the adjoint group representation:
\begin{align}
  \label{eq:retraction:tangent:adjoint}
  \diff\tau_\xi ={}& \Ad_{\tau(\xi)} \, \diff\tau_{-\xi} \,,&
  \diff\tau^{-1}_\xi ={}& \diff\tau^{-1}_{-\xi} \, \Ad_{\tau(-\xi)} \,.
\end{align}

\subsection{The group of rotations in \texorpdfstring{\(\bbR^3\)}{R³}}

The special orthogonal group of \(\bbR^3\), denoted \(SO(3)\), is the set of rotations of \(\bbR^3\) which can be identified with the group of orthogonal \(3\times3\) matrices with positive determinant. Other possible identifications are with the real projective space \(\mathbb{P}^{3}(\bbR)\), or with the closed ball of radius \(\pi\) whose surface is ``glued'' together at antipodal points. A vector in such set identifies with the axis of the rotation and its length gives the rotation angle, being 0 the identity.

The Lie algebra associated to \(SO(3)\) (and \(O(3)\)), denoted \(\frakso(3)\), consists (under identification) of the skew-symmetric \(3\times3\) matrices. Besides of the exponential map, which (under identification) corresponds here to the matrix exponential (\Cref{sec:expm}), another example of retraction is the Cayley transform (\Cref{sec:cayley}).


\section{Matrix identities}
\label{sec:matrix}
We summarize here some identities that relate common operations in \(\bbR^3\): the scalar product, the tensor product, the cross product, and the hat map. We recall that
\(
  x=(x_1,x_2,x_3)\in\bbR^3 \longmapsto \hat x =
  \left(
    \begin{smallmatrix}0&-x_3&x_2\\x_3&0&-x_1\\-x_2&x_1&0\end{smallmatrix}
  \right)
  \in\frakso(3)
\),
with inverse \((\hat x)^\vee=x\).

\begin{subequations}
  \label{eq:hats}
  \begin{gather}
    \allowdisplaybreaks
    \label{eq:xtensorxy}
    (x\otimes x)(y)=\langle x,y\rangle x\\
    \label{eq:xtensorx}
    x\otimes x = \|x\|^2I+\hat x^2\\
    \label{eq:xcrossy}
    \hat xy=x\times y\\
    \label{eq:hats:xy}
    \hat x\hat y=y\otimes x-\langle x,y\rangle I\\
    \label{eq:hats:xyyx}
    \hat x\hat y-\hat y\hat x=\widehat{x\times y}\\
    \label{eq:hats:xyx}
    \hat x\hat y\hat x=-\langle x,y \rangle\hat x\\
    \label{eq:hats:x2yyx2}
    \hat x^2\hat y+\hat y\hat x^2=-\|x\|^2\hat y-\langle x,y \rangle\hat x\\
    \label{eq:hats:x3}
    \hat x^3 = -\|x\|^2\hat x\\
    \label{eq:hats:trx2}
    \tr(\hat x^2) = -2\|x\|^2\\
    \label{eq:tr-skew}
    \big(\tfrac{\tr(R)-1}2\big)^2+\|(R^-)^\vee\|^2=1
  \end{gather}
\end{subequations}


\section{The Cayley transform}
\label{sec:cayley}
%
%
The Cayley transform is the map
\begin{equation}
  \label{eq:cay:def}
  \cay\colon \hat x\in\frakso(3) \longmapsto (I-\hat x)^{-1}(I+\hat x)\in SO(3) \,.
\end{equation}
Indeed, \(\cay(\hat x)^t\cay(\hat x)=I\).
We then have the formulas (see also \cite[Appendix B]{IserlesEtAl00}):
\begin{subequations}
  \label{eq:cay:formulas}
  \begin{align}
    \label{eq:cay}
    \cay(\hat x) ={}& I+2\lambda\hat x+2\lambda\hat x^2 \,,\\
    \label{eq:cay:inv}
    \cay^{-1}(R) ={}& \tfrac2{1+\tr(R)} R^- \,,\\
    \label{eq:cay:d}
    \dcay(\hat x) ={}& 2\lambda(I\pm\hat x) \,,\\
    \label{eq:cay:dinv}
    \dcay^{-1}(\hat x) ={}& \tfrac12(I\mp\hat x+x\otimes x) \,,
  \end{align}
\end{subequations}
where \(\lambda\coloneq\tfrac1{1+\|x\|^2}\), \(\dcay(\hat x)(y)\coloneq\left((\lra_{\cay(\hat x)^{-1}}\circ\Tan_{\hat x}\cay)(\hat y)\right)^\vee\), and \(\dcay^{-1}(\hat x)\coloneq(\dcay(\hat x))^{-1}\).
The lower\textbackslash{}upper signs in \cref{eq:cay:d,eq:cay:dinv} correspond to the choice \(\lra=\la\backslash{}\ra\), the left\textbackslash{}right action, respectively.
\[
  \xymatrix{
    \frakso(3)\cong\Tan_{\hat x}\frakso(3)\ar[r]^{\Tan_{\hat x}\cay} &
    \Tan_{\cay(\hat x)}SO(3) \ar[r]^{\lra_{\cay(\hat x)^{-1}}} &
    \Tan_ISO(3)\cong\frakso(3) \ar[d]^{\vee}\\
    \bbR^3 \ar@<+.7ex>[rr]^{\dcay(\hat x)} \ar[u]^{\wedge}&&
    \bbR^3 \ar@<+.3ex>[ll]^{\dcay^{-1}(\hat x)}
  }
\]

To prove the above formulas, we first show that
\begin{equation}
  \label{eq:invImX}
  (I-\hat x)^{-1} = I+\lambda\hat x+\lambda\hat x^2 \,.
\end{equation}
Indeed, carry out the following product and use \cref{eq:hats:x3} to get
\[
  (I+\lambda\hat x+\lambda\hat x^2)(I-\hat x)
  = I+\lambda\hat x+\lambda\hat x^2-\hat x-\lambda\hat x^2-\lambda\hat x^3
  = I+(\lambda-1+\lambda\|x\|^2)\hat x
  = I \,.
\]
An almost identical development using now \cref{eq:invImX} in definition \cref{eq:cay:def} gives formula \cref{eq:cay},
\[
  (I+\lambda\hat x+\lambda\hat x^2)(I+\hat x)
  = I+\lambda\hat x+\lambda\hat x^2+\hat x+\lambda\hat x^2+\lambda\hat x^3
  = I+(\lambda+1-\lambda\|x\|^2)\hat x + 2\lambda\hat x^2 \,.
\]
Besides, \cref{eq:invImX} also gives the commutativity of the factors in \cref{eq:cay:def},
\[
  (I-\hat x)^{-1}(I+\hat x) = (I+\hat x)(I-\hat x)^{-1} \,.
\]

For the inverse transform, \(\cay^{-1}\), define \(R \coloneq \cay(\hat x)\) to obtain thanks to Equations \cref{eq:cay,eq:hats:trx2}
\[
  \tr(R)
  = \tr(I)+2\lambda\tr(\hat x^2)
  = 3-4\lambda\|x\|^2
  = 4\lambda-1 \,,
\]
or, equivalently,
\[
  2\lambda = \tfrac12(1+\tr(R)) \,.
\]
Since \(R^- = 2\lambda\hat x\), we deduce formula \cref{eq:cay:inv}.
Moreover, from this same formula and the definition of \(\lambda\), we get the relation
\[
  1+\|x\|^2 \eqcolon \lambda^{-1} = 2\cdot\tfrac2{1+\tr(R)} \,.
\]
Taking into account that now we have \(\|x\| = \tfrac2{1+\tr(R)}\|(R^-)^\vee\|\) from \cref{eq:cay:inv}, we get
\[
  1+(\tfrac2{1+\tr(R)})^2\|(R^-)^\vee\|^2 = 2\cdot\tfrac2{1+\tr(R)} \,.
\]
Multiply by \((\tfrac{1+\tr(R)}2)^2\), pull everything to the left hand side,
\[
  \left(\tfrac{1+\tr(R)}2\right)^2 - 2\cdot\tfrac{1+\tr(R)}2 + \|(R^-)^\vee\|^2 = 0 \,,
\]
and complete squares to obtain the trigonometric relation \cref{eq:tr-skew} between the trace of a rotation and the norm of its skewsymmetric part.



For the tangent map, simple derivation yields
\begin{align*}
  (\Tan_{\hat x}\cay)(\hat y)
  \coloneq{}& \textstyle\ddt\left[\cay(\hat x(t))\right]|_{t=0} \ :\ \hat x(0)=\hat x\ \&\ \ddt\hat x(t)|_{t=0}=\hat y\\
  ={}& (I-\hat x)^{-1}\hat y(I-\hat x)^{-1}(I+\hat x) + (I-\hat x)^{-1}\hat y\\
  ={}& (I-\hat x)^{-1}\hat y(\cay(\hat x)+I)\\
  ={}& 2(I-\hat x)^{-1}\hat y(I-\hat x)^{-1} \,.
\end{align*}
Pulling to the identity by the left action (right action is analogous) results in
\begin{align*}
  \widehat{\dcay(\hat x)(y)}
  \coloneq{}& \cay(\hat x)^{-1}\cdot(\Tan_{\hat x}\cay)(\hat y)\\
  ={}& \cay(-\hat x)\cdot(\Tan_{\hat x}\cay)(\hat y)\\
  ={}& 2(I+\hat x)^{-1}\hat y(I-\hat x)^{-1} \,.
\end{align*}
Instead of developing this expression, we work around it by computing first its inverse,
\begin{align*}
  \widehat{\dcay^{-1}(\hat x)(y)}
  \coloneq{}& \widehat{\dcay(\hat x)^{-1}(y)}\\
  ={}& \tfrac12(I+\hat x)\hat y(I-\hat x)\\
  ={}& \tfrac12(\hat y+\hat x\hat y-\hat y\hat x-\hat x\hat y\hat x)\\
  ={}& \tfrac12(\hat y+\widehat{x\times y}+\langle x,y\rangle\hat x) \,.
\end{align*}
Equations \cref{eq:xcrossy,eq:xtensorxy} show the desired result, \cref{eq:cay:dinv},
which in turn is used in conjunction with \cref{eq:xtensorx,eq:invImX} to show \cref{eq:cay:d},
\begin{align*}
  \dcay^{-1}(\hat x)
  ={}& \tfrac12(I+\hat x+x\otimes x) \,\\
  ={}& \tfrac12((1+\|x\|^2)I+\hat x+\hat x^2)\\
  ={}& \tfrac1{2\lambda}(I-\hat x)^{-1} \,.
\end{align*}


\section{The matrix exponential in \texorpdfstring{\boldmath\(\frakso(3)\)}{so(3)}}
\label{sec:expm}
The matrix exponential is the map
\begin{equation}
  \label{eq:exp:def}
  \exp\colon A\in\frakgl(n) \longmapsto \sum_{k=0}^\infty\frac{A^k}{k!}\in GL(n) \,,
\end{equation}
whose restriction to \(\frakso(3)\) gives a map \(\exp\colon\frakso(3) \to SO(3)\).
%
We then have the formulas (see also \cite[Appendix B]{IserlesEtAl00})\footnote{Be aware of a typo in \cite[Eq.~(B.11)]{IserlesEtAl00}.}:
\begin{subequations}
  \label{eq:exp:formulas}
  \begin{align}
    \label{eq:exp}
    \exp(\hat x) ={}& I+\tfrac{\sin\omega}{\omega}\hat x+\tfrac{1-\cos\omega}{\omega^2}\hat x^2\\
    \label{eq:exp:inv}
    \log(R) ={}& \frac{\cos^{-1}\left(\tfrac{\tr(R)-1}2\right)}{\|(R^-)^\vee\|}R^-\\
    \label{eq:exp:d}
    \dexp(\hat x) ={}& I \pm \tfrac12\tfrac{\sin(\sfrac\omega2)}{(\sfrac\omega2)^2}\hat x + \tfrac{\omega-\sin(\omega)}{\omega^3}\hat x^2\\
    \label{eq:exp:dinv}
    \dlog(\hat x) ={}& I \mp \tfrac12\hat x+\tfrac12\tfrac{2-\omega\cot(\sfrac\omega2)}{\omega^2}\hat x^2
  \end{align}
\end{subequations}
where \(\omega=\|x\|\).
As in \cref{eq:cay:formulas}, the lower\textbackslash{}upper signs in \cref{eq:exp:d,eq:exp:dinv} correspond to the choice \(\lra=\la\backslash{}\ra\), the left\textbackslash{}right action, respectively.
\[
  \xymatrix{
    \frakso(3)\cong\Tan_{\hat x}\frakso(3)\ar[r]^{\Tan_{\hat x}\exp} &
    \Tan_{\exp(\hat x)}SO(3) \ar[r]^{\lra_{\exp(\hat x)^{-1}}} &
    \Tan_ISO(3)\cong\frakso(3) \ar[d]^{\vee}\\
    \bbR^3 \ar@<+.7ex>[rr]^{\dexp(\hat x)} \ar[u]^{\wedge}&&
    \bbR^3 \ar@<+.3ex>[ll]^{\dlog(\hat x)}
  }
\]

Formula \cref{eq:exp} is easily obtained by splitting the exponential series in odd and even terms so that the sine and cosine series are recovered.
\begin{align*}
  \exp(\hat x)
  ={}& \sum_{k=0}^\infty\frac{\hat x^k}{k!}\\
  ={}& I + \sum_{k=0}^\infty\frac{\hat x^{2k+1}}{(2k+1)!} + \sum_{k=0}^\infty\frac{\hat x^{2k+2}}{(2k+2)!}\\
  ={}& I + \sum_{k=0}^\infty(-1)^k\frac{\omega^{2k}}{(2k+1)!}\hat x + \sum_{k=0}^\infty(-1)^k\frac{\omega^{2k}}{(2k+2)!}\hat x^2\\
  ={}& I + \frac1\omega\left(\sum_{k=0}^\infty(-1)^k\frac{\omega^{2k+1}}{(2k+1)!}\right)\hat x - \frac1{\omega^2}\left(\sum_{k=0}^\infty(-1)^{k+1}\frac{\omega^{2k+2}}{(2k+2)!}\right)\hat x^2
\end{align*}

For the logarithm, define \(R \coloneq \exp(\hat x)\). Formula \cref{eq:exp} readily gives
\( R^- = \tfrac{\sin\omega}{\omega}\hat x \),
from which
\[
  \log(R) = \hat x = \tfrac{\omega}{\sin\omega}R^-
  \qquad\text{and}\qquad
  |\sin\omega| = \|(R^-)^\vee\|
  \,.
\]
Also from \cref{eq:exp}, and using \cref{eq:hats:trx2}, we get
\(
  \tr(R)
  = \tr(I)+\tfrac{1-\cos\omega}{\omega^2}\tr(\hat x^2)
  = 1+2\cos\omega
\)
or, equivalently,
\[
  \cos\omega = \tfrac{\tr(R)-1}2 \,,
\]
which show \cref{eq:exp:inv} for \(\omega\in[0,\pi]\). Besides, the trigonometric relations give \cref{eq:tr-skew} too.

For the time being, let
\[
  a(\omega)\coloneq\tfrac{\sin\omega}{\omega} \,,\qquad
  b(\omega)\coloneq\tfrac{1-\cos\omega}{\omega^2} \,,\quad\text{and}\quad
  \gamma\coloneq\langle x,y\rangle/\omega,
\]
so that simple derivation yields for the tangent map
\begin{align*}
  (\Tan_{\hat x}\exp)(\hat y)
  \coloneq{}& \textstyle\ddt\left[\exp(\hat x(t))\right]|_{t=0} \ :\ \hat x(0)=\hat x\ \&\ \ddt\hat x(t)|_{t=0}=\hat y\\
  ={}& a'\cdot\gamma\cdot\hat x + a\cdot\hat y + b'\cdot\gamma\cdot\hat x + b\cdot(\hat x\hat y+\hat y\hat x) \,.
\end{align*}
Pulling to the identity by the left action (right action is analogous) results in
\begin{align*}
  \widehat{\dexp(\hat x)(y)}
  \coloneq{}& \exp(\hat x)^{-1}\cdot(\Tan_{\hat x}\exp)(\hat y)\\
  ={}& \exp(-\hat x)\cdot(\Tan_{\hat x}\exp)(\hat y)\\
  ={}& a'\gamma\hat x + a\hat y + b'\gamma\hat x^2 + b\hat x\hat y + b\hat y\hat x\\
     & - aa'\gamma\hat x^2 - a^2\hat x\hat y - ab'\gamma\hat x^3 - ab\hat x^2\hat y - ab\hat x\hat y\hat x\\
     & + a'b\gamma\hat x^3 + ab\hat x^2\hat y + bb'\gamma\hat x^4 + b^2\hat x^3\hat y + b^2\hat x^2\hat y\hat x\\
  ={}& (a'+ab'\omega^2+ab\omega-a'b\omega^2)\gamma\hat x + a\hat y - \tfrac12(a^2+b^2\omega^2)(\hat x\hat y-\hat y\hat x)\\
  ={}& (a'/\omega + ab'\omega+ab-a'b\omega)\langle x,y\rangle\hat x + a\hat y - \tfrac12(a^2+b^2\omega^2)\widehat{x\times y}
\end{align*}
From formulas \cref{eq:xtensorxy,eq:xcrossy}, we may write
\begin{align*}
  \dexp(\hat x)(y)
  ={}& aI - \tfrac12(a^2+b^2\omega^2)\hat x  + (\tfrac{a'}\omega + ab'\omega+ab-a'b\omega)x\otimes x\\
\intertext{which is simplified using the expressions of \(a\) and \(b\) to get}
  ={}& aI - b\hat x  + \tfrac{1-a}{\omega^2}x\otimes x\\
  ={}&  I - b\hat x  + \tfrac{1-a}{\omega^2}\hat x^2
\end{align*}
where \cref{eq:xtensorx} has been used.

For its inverse \cref{eq:exp:dinv}, we take a direct approach by developing
\begin{align*}
  \left(I + \tfrac12\hat x + \tfrac1{\omega^2}(1-\tfrac12\tfrac{a}b)\hat x^2\right)\dexp(\hat x)
  ={}& I + \tfrac1{\omega^2}(1-\tfrac12b\omega^2 -\tfrac12\tfrac{a^2}b)\hat x^2 = I \,,
\end{align*}
where the last term cancels proving the desired result.


\section{The skewsymmetric matrix projection}
\label{sec:skew}
The skewsymmetric matrix projection is the linear endomorphism
\begin{equation}
  \label{eq:skew:def}
  \skew\colon A\in\mathcal{M}_{n\times n}(\bbR) \longmapsto A^-\coloneq\tfrac12(A-A^t)\in\mathcal{M}_{n\times n}(\bbR) \,.
\end{equation}
This map is indeed a projection that annihilates symmetric matrices and, therefore, it is not bijective. Its restriction to \(SO(n)\) is however a local diffeomorphism around the identity whose inverse is the retraction map
\begin{equation}
  \label{eq:unskew:def}
  \unskew\colon A\in\frakso(n) \longmapsto A + \sqrt{I+A^2} \in SO(n) \,,
\end{equation}
where \(\sqrt{I+A^2}\) is the unique positive definite matrix whose square is \(I+A^2\) for \(A\) small enough \cite[Thm.~6.1]{CardosoLeite03}.
For the case \(n=3\), we have the following formulas
\begin{subequations}
  \label{eq:skew:formulas}
  \begin{align}
    \label{eq:unskew}
    \unskew(\hat x) ={}& I + \hat x + \gamma\hat x^2\\
    \label{eq:skew}
    \skew(R) ={}& \tfrac12(R-R^t)\\
    \label{eq:unskew:d}
    \dunskew(\hat x) ={}& \gamma I \pm \tfrac\gamma{3+2\gamma}\hat x + \tfrac{1+\gamma^2}{3+2\gamma}\hat x^2\\
    \label{eq:skew:d}
    \dskew(\hat x) ={}& \gamma^{-1} I \mp \tfrac12\hat x - \tfrac12\gamma\hat x^2
  \end{align}
\end{subequations}
where \(\gamma^{-1}=1+\sqrt{1-\|x\|^2}\) (so here ``small enough'' means \(0\leq\|x\|<1\)).
As in \cref{eq:cay:formulas,eq:exp:formulas}, the lower\textbackslash{}upper signs in \cref{eq:unskew:d,eq:skew:d} correspond to the choice \(\lra=\la\backslash{}\ra\), the left\textbackslash{}right action, respectively.
\[
  \xymatrix{
    \frakso(3)\cong\Tan_{\hat x}\frakso(3) &
    \Tan_{R}SO(3) \ar[l]_{\quad\ \Tan_R\skew} &
    \Tan_ISO(3)\cong\frakso(3) \ar[l]_{\lra_R\quad} \ar[d]^{\vee}\\
    \bbR^3 \ar@<+.7ex>[rr]^{\dunskew(\hat x)} \ar[u]^{\wedge} &&
    \bbR^3 \ar@<+.3ex>[ll]^{\dskew(\hat x)}
  }
\]

Equation \cref{eq:unskew} is easily proven if we observe that, for \(x\in\bbR^3\) small enough, \(\sqrt{I+\hat x^2}=I+\gamma\hat x^2\) since, by Eq. \cref{eq:hats:x3}, \((I+\gamma\hat x^2)^2 = I+\hat x^2 \) and \(I+\gamma\hat x^2\) is positive definite.
Next we show Equation \cref{eq:skew:d}. To this end, take \(\hat x = \skew(R)\) and compute
\begin{align*}
  \widehat{\dskew(\hat x)(y)}
  ={}& \Tan_R\skew\big(\hat y\unskew(\hat x)\big)\\
  ={}& \skew\big(\hat y\,(I + \hat x + \gamma\hat x^2)\big)\\
  ={}& \hat y + \tfrac12(\hat y\hat x-\hat x\hat y) + \tfrac12\gamma(\hat y\hat x^2+\hat x^2\hat y)\\
  ={}& \hat y - \tfrac12\widehat{x\times y} - \tfrac12\gamma(\|x\|^2\hat y+\langle x,y \rangle\hat x)
\end{align*}
which shows
\[ \dskew(\hat x) = I - \tfrac12 \hat x - \tfrac12\gamma(\|x\|^2I+x\otimes x) \,. \]
In this, we have used the fact that \(\Tan_R\skew=\skew\) and the matrix identities %
\cref{eq:hats}, which in turn give \cref{eq:skew:d}.

To show that \cref{eq:skew:d} is the inverse of \cref{eq:unskew:d}, simply expand the matrix product of both to get the identity.


\section{Continuous and discrete Euler-Lagrange equations for Lie groups}
\label{sec:euler-lagrange}
In this section we recall the continuous and discrete Euler-Lagrange equations for systems with configuration a Lie group \(G\) (with Lie algebra \(\frakg\)). In the continuous case the equations are determined prescribing a Lagrangian function \(L\colon\bbR\times TG\to\bbR\) and, in the discrete case, by a discrete Lagrangian function \(l\colon\bbZ\times G\times G\to\bbR\).

\subsection{The continuous equations}
Given a smooth manifold \(Q\), let \((q^i,\dot{q}^i)\) denote adapted coordinates on its tangent bundle \(TQ\), a Lagrangian function \(L\colon\bbR\times TQ\to\bbR\), and  an external force \(F\colon\bbR\times TQ\to T^*Q\) (a fibered map over \(Q\)). The Euler-Lagrange equations for the system \((L,F)\) are
\begin{equation}\label{eq:EL:group:continuous}
  \ddt\left(\pp[L]{\dot{q}}\right) - \pp[L]q = F \,.
\end{equation}
These equations are still valid for a Lie group \(G\), however they are usually rewritten in terms of left or right trivialization, \(TG\cong G\times\frakg\).

Given \(L\colon\bbR\times TG\to\bbR\) and \(F\colon\bbR\times TG\to T^*G\), define their right-trivializations \(\bar{L}\colon\bbR\times G\times\frakg\to\bbR\) and \(\bar{F}\colon\bbR\times G\times\frakg\to\frakg^*\) by the expressions
\begin{align*}
  \bar{L}(t,g,\xi) ={}& L(t,g,\ra_g\xi) \,,
  & \bar{F}(t,g,\xi) ={}& \ra_g^*\big(F(t,g,\ra_g\xi)\big) \,.
\end{align*}
The right-trivialized Euler-Lagrange equation for \((\bar{L},\bar{F})\) are
\begin{equation}\label{eq:EL:group:continuous:reduced}
  \ddt\left(\pp[\bar{L}]\xi\right)
  + \ad_{\xi}^*\left(\pp[\bar{L}]\xi\right)
  - {\ra}_g^*\left(\pp[\bar{L}]g\right) = \bar{F} \,,
\end{equation}
which, together with the \emph{reconstruction equation} \(\dot{g}=\xi g\), are equivalent to Eq. \cref{eq:EL:group:continuous}.

Given smooth functions \(a,b\colon\bbR\to\bbR_+\) and \(\phi\colon G\to\bbR\), consider the (right) trivialized Lagrangian
\[
  \bar{L}(t,g,\xi) = \tfrac{1}{2}a(t)\norm{\xi}^2 - b(t)\phi(g) \,,
\]
where \(\norm{\!\cdot\!}\) is the norm associated to a given inner product \(\langle\,\cdot\!\;,\cdot\!\;\rangle\) on the Lie algebra \(\frakg\).
Then the right-trivialized Euler-Lagrange equation are
\[
  \ddt\left(a(t)\xi^\flat\right)
  + a(t)\ad_{\xi}^*\xi^\flat - b(t)\ra_g^*\diff\phi(g) = 0 \in\frakg^* \,,
\]
that is, after expanding and using the sharp isomorphism,
\[
  \dot{\xi}+\ad_{\xi}^t\xi
  + \tfrac{\dot{a}}{a}\,\xi  - \tfrac{b}{a}\nabla\phi(g) = 0 \in\frakg \,.
\]

\subsection{The discrete equations}
In this case, the phase space \(TQ\) is replaced by \(Q\times Q\), while the continuous time line \(\bbR\) is replaced by discrete time ticks \(\bbZ\). We therefore consider a time-dependent discrete Lagrangian \(l\colon\bbZ\times Q\times Q\to\bbR\), otherwise a family \(l_k\colon Q\times Q\to\bbR\), \(k\in\bbZ\), and two families of external forces \(f^\pm_k\colon Q\times Q\to T^*Q\) (fibered maps over \(Q\) along the projections \(\text{pr}_{\pm}\)). Then the discrete Euler-Lagrange equations for the system \((l,f^\pm)\) are (confer with \cite{CMM23}, for this approach, and with \cite{marsden-west}, for an introduction to discrete Lagrangian mechanics):
\begin{equation}\label{eq:DEL:point-point}
  D_1l_k(q_k,q_{k+1}) + D_2l_{k-1}(q_{k-1},q_k) + f^-_k(q_k,q_{k+1}) + f^+_{k-1}(q_{k-1},q_k) = 0 \in \Tan_{q_k}^*Q\,.
\end{equation}
In this picture, given two initial points \((q_0,q_1)\), Eq. \cref{eq:DEL:point-point} determines iteratively \(q_{k+1}\) from the two previous points \((q_{k-1},q_k)\) for \(k\geq1\).

For the case where \(Q\) is a Lie group \(G\), in the spirit of the earlier trivialized expressions, instead of working with pairs \((g_k,g_{k+1})\) of consecutive points in a trajectory, one can chose to work with ``pointing arrows'', pairs of the form \emph{source-arrow} \((g_k,h_k)\) pointing towards a \emph{target} \(g_kh_k=g_{k+1}\). With this perspective in mind, define the ``trivialized'' discrete Lagrangian and forces as follows
\begin{align*}
  \bar{l}_k(g,h) \coloneq{}& l_k(g,gh) \,,
  & \bar{f}^\pm_k(g,h) \coloneq{}& \la_{\text{pr}_\pm(g,gh)}^*f^\pm_k(g,gh) \,,
\end{align*}
where \(\text{pr}_-\) and \(\text{pr}_+\) are the source and target projection, respectively.
After simple manipulation, together with the reconstruction equation \(g_{k+1}=g_kh_k\), the Euler-Lagrange equation \cref{eq:DEL:point-point} reads
\begin{equation}\label{eq:DEL:point-trajectory}
  \textstyle
  \la_{g_k}^*\partial_g\bar{l}_k - \ra_{h_k}^*\partial_h\bar{l}_k + \la_{h_{k-1}}^*\partial_h\bar{l}_{k-1} + \bar f^-_k + \bar f^+_{k-1} = 0 \in\frakg^*\,,
\end{equation}
where  \(\partial_{m}\bar{l}_k\) is a shorthand notation for \(\partial_m\bar{l}_k(g_k,h_k)\) with \(m=g,h\), and similarly for \(\bar{f}^\pm_k\).


\bibliographystyle{siamplain}
\bibliography{references}

\end{document}